\documentstyle[12pt]{article}    
\input{amssym.def}
\input{amssym.tex}
\begin{document}      
\title{Operator algebras and conjugacy problem for the 
pseudo-Anosov automorphisms of a surface}    

\author{Igor  ~Nikolaev
\footnote{Partially supported 
by NSERC.}}


\date{}    
 \maketitle    
    

\newtheorem{thm}{Theorem}    
\newtheorem{lem}{Lemma}    
\newtheorem{dfn}{Definition}    
\newtheorem{rmk}{Remark}    
\newtheorem{cor}{Corollary}    
\newtheorem{prp}{Proposition}    
\newtheorem{exm}{Example}    
\newtheorem{con}{Conjecture}
\newtheorem{prb}{Problem}


\newcommand{\tr}{\hbox{\bf Tr}}
\newcommand{\dete}{\hbox{\bf Det}}
\newcommand{\Sing}{\hbox{\bf Sing}}


\centerline{{\it  In memory of W.~P.~Thurston}}

 \begin{abstract}    

The conjugacy problem for the pseudo-Anosov automorphisms
of a compact surface is studied. To each pseudo-Anosov automorphism
$\phi$, we assign an AF $C^*$-algebra ${\Bbb A}_{\phi}$ (an operator 
algebra). It is proved that the assignment is functorial, i.e.
every $\phi'$, conjugate to $\phi$, maps to an AF $C^*$-algebra 
${\Bbb A}_{\phi'}$, which is stably isomorphic to ${\Bbb A}_{\phi}$. 
The new invariants of the conjugacy of the pseudo-Anosov automorphisms
are obtained from the known invariants of the stable isomorphisms
of the AF $C^*$-algebras. Namely, the main invariant is a triple
$(\Lambda, [I], K)$, where $\Lambda$ is an order in the ring of
integers in a real algebraic number field $K$ and $[I]$ an
equivalence class of the ideals in $\Lambda$. The numerical
invariants include the determinant $\Delta$ and the signature
$\Sigma$, which we compute for the case of the Anosov automorphisms. 
A question concerning the $p$-adic invariants of the pseudo-Anosov
automorphism is formulated.

\vspace{7mm}    
    

\vspace{5mm}
\end{abstract}

\section*{Introduction} 
{\bf A. ~Conjugacy problem.}
Let $Mod~(X)$ be the mapping class group of a compact  surface $X$, i.e. the 
group of orientation preserving automorphisms of $X$ modulo the trivial ones.
Recall that $\phi,\phi'\in Mod~(X)$ are  conjugate automorphisms,  whenever $\phi'=h\circ\phi\circ h^{-1}$
for an  $h\in Mod~(X)$. It is not hard to see  that  conjugation is an equivalence relation 
which splits the mapping class group  into  disjoint classes of conjugate automorphisms.  
The construction of invariants of the conjugacy classes in $Mod~(X)$ is an important and difficult 
problem studied by Hemion \cite{Hem1}, Mosher \cite{Mos1},  and others. 
Any knowledge of  such invariants leads to  a topological classification of  three-dimensional 
manifolds, which fiber over the circle with  monodromy $\phi\in Mod~(X)$ \cite{Thu2}.

\medskip\noindent
{\bf B. ~Pseudo-Anosov automorphisms.}
It is known that any $\phi\in Mod ~(X)$ is isotopic to an automorphism
$\phi'$, such that either (i) $\phi'$ has a finite order, or
(ii) $\phi'$ is a pseudo-Anosov (aperiodic) automorphism, or else
(iii) $\phi'$ is reducible by a system of curves $\Gamma$ surrounded
 by the small tubular neighborhoods $N(\Gamma)$, such that on
 $X \setminus  N(\Gamma)$ $\phi'$ satisfies either (i) or (ii).
Let $\phi$ be a representative of the equivalence class 
of a pseudo-Anosov automorphism. 
Then there exist a pair consisting of the stable ${\cal F}_s$
and unstable ${\cal F}_u$ mutually orthogonal measured foliations on the surface $X$,
such that $\phi({\cal F}_s)={1\over\lambda_{\phi}}{\cal F}_s$ 
and $\phi({\cal F}_u)=\lambda_{\phi}{\cal F}_u$, where $\lambda_{\phi}>1$
is called a dilatation of $\phi$. The foliations ${\cal F}_s,{\cal F}_u$ are minimal,
uniquely ergodic and describe the automorphism $\phi$ up to a power.  
In the sequel, we shall focus on the conjugacy problem for the pseudo-Anosov 
automorphisms of a surface $X$.

\medskip\noindent
{\bf C.  AF $C^*$-algebras.}
The $C^*$-algebra is an algebra $A$ over ${\Bbb C}$ with a norm
$a\mapsto ||a||$ and an involution $a\mapsto a^*$ such that
it is complete with respect to the norm and $||ab||\le ||a||~ ||b||$
and $||a^*a||=||a^2||$ for all $a,b\in A$. The $C^*$-algebras have been introduced by
Murray and von Neumann as rings of bounded operators on a 
Hilbert space and are strongly connected  with the geometry and topology
of manifolds \cite{B}, \S 24.  
Any  simple finite-dimensional $C^*$-algebra is isomorphic 
to the algebra $M_n({\Bbb C})$ of the complex $n\times n$ matrices. 
A natural completion of the finite-dimensional semisimple  $C^*$-algebras 
(as $n\to\infty$) is  known as an  {\it AF $C^*$-algebra} \cite{E}.  
AF $C^*$-algebra is most conveniently  given by an infinite graph,
which records the inclusion of the finite-dimensional subalgebras 
into the AF $C^*$-algebra. The graph is called  a {\it Bratteli diagram}.
When the diagram is periodic,  the  AF $C^*$-algebra is {\it stationary};
this is an important special case.  In  addition to the usual 
isomorphism $\cong$, the $C^*$-algebras $A, A'$
are called {\it stably isomorphic} whenever $A\otimes {\cal K}\cong A'\otimes {\cal K}$,
where ${\cal K}$ is the $C^*$-algebra of compact operators.

\medskip\noindent
{\bf D. ~Motivation.}
Let $\phi\in Mod~(X)$ be a pseudo-Anosov automorphism. The main idea
of present paper is to assign to $\phi$ an AF $C^*$-algebra, ${\Bbb A}_{\phi}$,
so that for  every $h\in Mod~(X)$  the following diagram commutes:

\begin{picture}(300,110)(-80,-5)
\put(20,70){\vector(0,-1){35}}
\put(130,70){\vector(0,-1){35}}
\put(45,23){\vector(1,0){53}}
\put(45,83){\vector(1,0){53}}
\put(15,20){${\Bbb A}_{\phi}$}
\put(128,20){${\Bbb A}_{\phi'}$}
\put(17,80){$\phi$}
\put(117,80){$\phi'=h\circ\phi\circ h^{-1}$}
\put(60,30){\sf stable}
\put(50,10){\sf isomorphism}
\put(54,90){\sf conjugacy}
\end{picture}

\noindent
(In other words, if $\phi$ and  $\phi'$ are  conjugate pseudo-Anosov automorphisms,  then
the AF $C^*$-algebras ${\Bbb A}_{\phi}$ and  ${\Bbb A}_{\phi'}$ are stably isomorphic.)  For the sake of 
clarity, we  shall consider  an example illustrating  the idea in the  case  $X=T^2$ (a torus).

\medskip\noindent
{\bf E.  ~Model example.}
Let $\phi\in Mod~(T^2)$ be the Anosov automorphism given by a non-negative matrix 
$A_{\phi}\in SL_2({\Bbb Z})$. (The assumption is not restrictive; each $A_{\phi}$
with $\tr~(A_{\phi})>0$ is similar to a non-negative matrix. The case $\tr~(A_{\phi})<0$
is treated likewise -- by reduction to a non-positive matrix; then the absolute value of 
all entries must be  taken.)  Consider a stationary AF $C^*$-algebra, 
${\Bbb A}_{\phi}$, given by the following periodic Bratteli diagram:

\begin{figure}[here]
\begin{picture}(350,100)(60,0)
\put(97,47){$\bullet$}

\put(117,27){$\bullet$}
\put(117,67){$\bullet$}

\put(157,27){$\bullet$}
\put(157,67){$\bullet$}

\put(197,27){$\bullet$}
\put(197,67){$\bullet$}

\put(237,27){$\bullet$}
\put(237,67){$\bullet$}


\put(100,50){\line(1,1){20}}
\put(100,50){\line(1,-1){20}}

\put(120,30){\line(1,0){40}}

\put(120,32){\line(1,1){40}}
\put(120,28){\line(1,1){40}}

\put(120,72){\line(1,0){40}}
\put(120,70){\line(1,0){40}}
\put(120,68){\line(1,0){40}}

\put(120,70){\line(1,-1){40}}


\put(160,30){\line(1,0){40}}

\put(160,32){\line(1,1){40}}
\put(160,28){\line(1,1){40}}

\put(160,72){\line(1,0){40}}
\put(160,70){\line(1,0){40}}
\put(160,68){\line(1,0){40}}

\put(160,70){\line(1,-1){40}}


\put(200,30){\line(1,0){40}}

\put(200,32){\line(1,1){40}}
\put(200,28){\line(1,1){40}}

\put(200,72){\line(1,0){40}}
\put(200,70){\line(1,0){40}}
\put(200,68){\line(1,0){40}}

\put(200,70){\line(1,-1){40}}

\put(250,30){$\dots$}
\put(250,70){$\dots$}

\put(137,78){$a_{11}$}
\put(177,78){$a_{11}$}
\put(217,78){$a_{11}$}


\put(116,55){$a_{12}$}
\put(156,55){$a_{12}$}
\put(196,55){$a_{12}$}


\put(113,44){$a_{21}$}
\put(153,44){$a_{21}$}
\put(193,44){$a_{21}$}


\put(137,22){$a_{22}$}
\put(177,22){$a_{22}$}
\put(217,22){$a_{22}$}

\put(290,50){
$A_{\phi}=\left(\matrix{a_{11} & a_{12}\cr a_{21} & a_{22}}\right)$,}

\end{picture}

\caption{The AF $C^*$-algebra  ${\Bbb A}_{\phi}$.}
\end{figure}
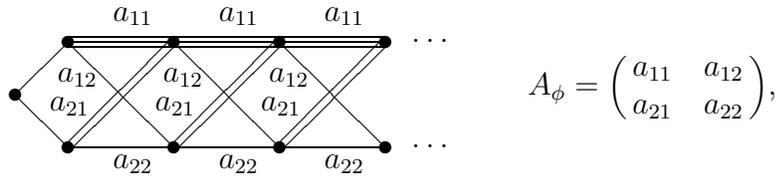

\noindent
where $a_{ij}$ indicate the multiplicity of the respective edges of the graph. 
We encourage  the reader to verify that $F: \phi\mapsto {\Bbb A}_{\phi}$
is a well-defined function on the set of Anosov automorphisms given by
the hyperbolic matrices with the non-negative entries.
Let us show that if $\phi,\phi'\in Mod~(T^2)$ are  conjugate Anosov automorphisms,
then ${\Bbb A}_{\phi},{\Bbb A}_{\phi'}$ are  stably isomorphic
AF $C^*$-algebras.  Indeed, let $\phi'=h\circ\phi\circ h^{-1}$ for
an $h\in Mod~(X)$. Then $A_{\phi'}=TA_{\phi}T^{-1}$ for
a matrix $T\in SL_2({\Bbb Z})$. Note that $(A_{\phi}')^n=(TA_{\phi}T^{-1})^n=
TA_{\phi}^nT^{-1}$, where $n\in {\Bbb N}$. We shall use the following
criterion:  the AF $C^*$-algebras ${\Bbb A},{\Bbb A}'$
are stably isomorphic if and only if their Bratteli diagrams contain a 
common block of an arbitrary length (compare with \cite{E}, Theorem 2.3; 
recall  that an order-isomorphism mentioned in the theorem is equivalent
to the condition that the corresponding Bratteli diagrams have the same
infinite tails -- i.e. a common block of infinite length). 
Consider the following sequences  of matrices:
$$
\left\{
\begin{array}{c}
\underbrace{A_{\phi}A_{\phi}\dots A_{\phi}}_n\\
T\underbrace{A_{\phi}A_{\phi}\dots A_{\phi}}_nT^{-1},
\end{array}
\right.
$$
which mimic the Bratteli diagrams of ${\Bbb A}_{\phi}$ and ${\Bbb A}_{\phi'}$.
Letting  $n\to\infty$, we conclude that  
${\Bbb A}_{\phi}\otimes {\cal K}\cong {\Bbb A}_{\phi'}\otimes {\cal K}$.

\medskip\noindent
{\bf F.  ~Invariants of torus  automorphisms obtained from the operator algebras.}
Conjugacy problem for the Anosov automorphisms can now be recast in terms of AF $C^*$-algebras: 
find  invariants of  stable isomorphism classes  of the stationary
AF $C^*$-algebras.  One such invariant is due to Handelman \cite{Han1}. 
Consider an eigenvalue problem for the hyperbolic matrix $A_{\phi}\in SL_2({\Bbb Z})$:
$A_{\phi}v_A=\lambda_Av_A$,  where $\lambda_A>1$ is the Perron-Frobenius eigenvalue  and 
$v_A=(v_A^{(1)},v_A^{(2)})$ the corresponding eigenvector with the positive entries 
normalized so that $v_A^{(i)}\in K={\Bbb Q}(\lambda_A)$. 
Denote by ${\goth m}={\Bbb Z}v_A^{(1)}+{\Bbb Z}v_A^{(2)}$
 the  ${\Bbb Z}$-module in the number field $K$. Recall that the coefficient
ring, $\Lambda$, of module ${\goth m}$ consists of the elements $\alpha\in K$
such that $\alpha {\goth m}\subseteq {\goth m}$. It is known that 
$\Lambda$ is an order in $K$ (i.e. a subring of $K$
containing $1$) and, with no restriction, one can assume that 
${\goth m}\subseteq\Lambda$. It follows  from the definition, that ${\goth m}$
coincides with an ideal, $I$, whose equivalence class in $\Lambda$ we shall denote
by $[I]$. It has been proved by Handelman, that the triple $(\Lambda, [I], K)$ is an arithmetic invariant of the 
stable isomorphism class of ${\Bbb A}_{\phi}$: the ${\Bbb A}_{\phi},{\Bbb A}_{\phi'}$
are stably isomorphic AF $C^*$-algebras if and only if $\Lambda=\Lambda', [I]=[I']$ and $K=K'$. 
It is interesting to compare the operator algebra invariants with the matrix invariants 
obtained in \cite{LaMa1} and \cite{Wal1}.

\medskip\noindent
{\bf G.  ~AF $C^*$-algebra ${\Bbb A}_{\phi}$ (pseudo-Anosov case).}
Denote by ${\cal F}_{\phi}$ the stable foliation of a pseudo-Anosov automorphism
$\phi\in Mod~(X)$. For brevity, we assume that ${\cal F}_{\phi}$ is an oriented 
foliation given by the trajectories of a closed $1$-form $\omega\in H^1(X; {\Bbb R})$.
Let $v^{(i)}=\int_{\gamma_i}\omega$, where $\{\gamma_1,\dots,\gamma_n\}$ is a basis
in the relative homology $H_1(X, \Sing~{\cal F}_{\phi}; {\Bbb Z})$, such that
$\theta=(\theta_1,\dots,\theta_{n-1})$ is a vector with  positive coordinates 
$\theta_i=v^{(i+1)} / v^{(1)}$. (Note that the $\theta_i$ depend on a basis in the 
homology group; but a ${\Bbb Z}$-module generated by the  $\theta_i$ does not --
see lemma \ref{lm1}.) 
Consider the (infinite) Jacobi-Perron continued
fraction \cite{BE} of $\theta$:
$$
\left(\matrix{1\cr \theta}\right)=
\lim_{k\to\infty} \left(\matrix{0 & 1\cr I & b_1}\right)\dots
\left(\matrix{0 & 1\cr I & b_k}\right)
\left(\matrix{0\cr {\Bbb I}}\right),
$$
where $b_i=(b^{(i)}_1,\dots, b^{(i)}_{n-1})^T$ is a vector of the nonnegative integers,  
$I$ the unit matrix and ${\Bbb I}=(0,\dots, 0, 1)^T$.
By the definition, ${\Bbb A}_{\phi}$ is an (isomorphism class of the)  AF $C^*$-algebra given by the 
Bratteli diagram   whose incidence matrices coincide with
$B_k=\left(\small\matrix{0 & 1\cr I & b_k}\right)$ for  $k=1,\dots, \infty$.
Note that this yields the Bratteli diagram derived in the model example (the Anosov case).

\medskip\noindent
{\bf H. ~Main results.}
For a matrix $A\in GL_n({\Bbb Z})$ with positive
entries, we denote by $\lambda_A$  the Perron-Frobenius eigenvalue
and let $(v^{(1)}_A,\dots, v^{(n)}_A)$ denote the corresponding normalized eigenvector
with  $v^{(i)}_A\in K={\Bbb Q}(\lambda_A)$.
The coefficient (endomorphism) ring of the module   ${\goth m}={\Bbb Z}v^{(1)}_A+\dots+{\Bbb Z}v^{(n)}_A$ 
will be denoted by $\Lambda$.  The equivalence class of ideal $I$ in $\Lambda$
 will be denoted  $[I]$. 
Finally, we denote by $\Delta=\dete~(a_{ij})$ and $\Sigma$ the determinant
and signature of the symmetric bilinear form $q(x,y)=\sum_{i,j}^na_{ij}x_ix_j$,
where $a_{ij}=\tr~(v^{(i)}_Av^{(j)}_A)$ and $\tr~ (\bullet)$ 
the trace function.  Our main results can be expressed as follows. 
\begin{thm}\label{thm1}
${\Bbb A}_{\phi}$ is a stationary AF $C^*$-algebra. 
\end{thm}
Let $\Phi$ be a category of all pseudo-Anosov (Anosov, resp.) automorphisms
of a surface of the genus $g\ge 2$ ($g=1$,  resp.);  the arrows (morphisms)
are conjugations between the automorphisms. Likewise,
let ${\cal A}$ be the category of all  stationary  AF $C^*$-algebras ${\Bbb A}_{\phi}$,
where  $\phi$ runs over the set $\Phi$; the arrows of ${\cal A}$ are  stable 
isomorphisms among the algebras ${\Bbb A}_{\phi}$. 
\begin{thm}\label{thm2}
Let  $F:\Phi\to {\cal A}$  be a map given by the formula $\phi\mapsto {\Bbb A}_{\phi}$.   Then:

\medskip
(i) $F$ is a functor;   it   maps  conjugate pseudo-Anosov automorphisms to  stably
isomorphic AF $C^*$-algebras;

\smallskip
(ii) $Ker~F=[\phi]$, where $[\phi]=\{\phi'\in \Phi~|~(\phi')^m=\phi^n, ~m,n\in {\Bbb N}\}$
is the commensurability class of the pseudo-Anoov automorphism $\phi$.  
\end{thm}
\begin{cor}\label{cor1}
The following are invariants of the conjugacy classes of the pseudo-Anosov 
automorphisms:

\medskip
(i) triples  $(\Lambda, [I], K)$;

\smallskip
(ii) integers $\Delta$ and $\Sigma$.  
\end{cor}

\medskip\noindent
{\bf I.  ~How to calculate invariants  $(\Lambda, [I], K)$, $\Delta$ and $\Sigma$?}
There is no easy way; the problem is comparable to that of numerical invariants 
of the fundamental group of a knot. A step in this direction would be 
computation of the matrix $A$; the latter is similar  to the matrix $\rho(\phi)$,
where $\rho: Mod~(X)\to PIL$ is a faithful representation of the mapping class
group as a group of the piecewise-integral-linear (PIL) transformations
\cite{Pen1}, p.45. The entries of $\rho(\phi)$ are the linear combinations 
of the Dehn twists along the $(3g-1)$ (Lickorish) curves on the surface $X$. 
Then one can effectively determine whether the $\rho(\phi)$ and $A$ are
similar matrices (over ${\Bbb Z}$) by bringing the polynomial matrices
$\rho(\phi)-xI$ and $A-xI$ to the Smith normal form; when the similarity
is established, the numerical invariants $\Delta$ and $\Sigma$ become the
polynomials in the Dehn twists. A tabulation of the simplest elements 
of $Mod~(X)$ is possible in terms of $\Delta$ and $\Sigma$ (see \S 4.3);
however, this task lies beyond the scope of present paper.

\medskip\noindent
{\bf K. ~Structure of the paper.}
Proofs of the main results can be found in section 3. 
Sections 1 and 2 consist of lemmas used to prove the
main results. Section 4 includes some examples, open
problems and conjectures. 
Since the paper does not include a formal section
on the preliminaries, we encourage the reader to consult
\cite{B}, \cite{E}, \cite{Kri1} (operator algebras \& dynamics),
\cite{HuMa1}, \cite{Thu1} (measured foliations) and 
\cite{BE}, \cite{Per1} (Jacobi-Perron continued fractions).

\tableofcontents

\section{Jacobian of a measured foliation}
\subsection{Definition of the jacobian}
Let ${\cal F}$ be a measured foliation on a compact surface
$X$ \cite{Thu1}. For the sake of brevity, we shall always assume
that ${\cal F}$ is an oriented foliation, i.e. given by the 
trajectories of a closed $1$-form $\omega$ on $X$. (The assumption
is not a  restriction -- each measured foliation is oriented on a 
surface $\widetilde X$, which is a double cover of $X$ ramified 
at the singular points of the half-integer index of the non-oriented
foliation \cite{HuMa1}.) Let $\{\gamma_1,\dots,\gamma_n\}$ be a
basis in the relative homology group $H_1(X, \Sing~{\cal F}; {\Bbb Z})$,
where $\Sing~{\cal F}$ is the set of singular points of the foliation ${\cal F}$.
It is well known that $n=2g+m-1$, where $g$ is the genus of $X$ and $m= |\Sing~({\cal F})|$. 
The periods of $\omega$ in the above basis will be written
$$
\lambda_i=\int_{\gamma_i}\omega.
$$
The real numbers $\lambda_i$ are coordinates of ${\cal F}$ in the space of all 
measured foliations on $X$ (with the fixed set of  singular points)
\cite{DoHu1}.  
\begin{dfn}
By a jacobian $Jac~({\cal F})$ of the measured foliation ${\cal F}$, we
understand a ${\Bbb Z}$-module ${\goth m}={\Bbb Z}\lambda_1+\dots+{\Bbb Z}\lambda_n$
regarded as a subset of the real line ${\Bbb R}$. 
\end{dfn}
An importance of the jacobians stems from an observation that
although the periods, $\lambda_i$,  depend on the choice of basis in 
$H_1(X, \Sing~{\cal F}; {\Bbb Z})$,  the jacobian does not.
Moreover,  up to a scalar multiple,   the jacobian is an invariant
of the equivalence class of the foliation ${\cal F}$. 
We formalize these observations in the following two sections.

\subsection{Invariance of  the jacobian}
\begin{lem}\label{lm1}
The ${\Bbb Z}$-module ${\goth m}$  is independent of choice of 
basis
\linebreak
in $H_1(X, \Sing~{\cal F}; {\Bbb Z})$ and depends solely on the foliation ${\cal F}$.
\end{lem}
{\it Proof.} Indeed, let $A=(a_{ij})\in GL_n({\Bbb Z})$ and let
$$
\gamma_i'=\sum_{j=1}^na_{ij}\gamma_j
$$
be a new basis in $H_1(X, \Sing~{\cal F}; {\Bbb Z})$. 
Then using the integration rules:
\begin{eqnarray}
\lambda_i'  &= \int_{\gamma_i'}\omega &= \int_{\sum_{j=1}^na_{ij}\gamma_j}\omega=\nonumber \\
 &= \sum_{j=1}^n\int_{\gamma_j}\omega  &= \sum_{j=1}^na_{ij}\lambda_j.\nonumber 
\end{eqnarray}

\bigskip
To prove that ${\goth m}={\goth m}'$, consider the following equations:
\begin{eqnarray}
{\goth m}'  &= \sum_{i=1}^n{\Bbb Z}\lambda_i' &= \sum_{i=1}^n {\Bbb Z} \sum_{j=1}^n a_{ij}\lambda_j=\nonumber \\
 &= \sum_{j=1}^n \left(\sum_{i=1}^n a_{ij}{\Bbb Z}\right)\lambda_j  &\subseteq  {\goth m}. \nonumber
\end{eqnarray}
Let $A^{-1}=(b_{ij})\in GL_n({\Bbb Z})$ be an inverse to the matrix $A$.
Then $\lambda_i=\sum_{j=1}^nb_{ij}\lambda_j'$ and 
\begin{eqnarray}
{\goth m}  &= \sum_{i=1}^n{\Bbb Z}\lambda_i &= \sum_{i=1}^n {\Bbb Z} \sum_{j=1}^n b_{ij}\lambda_j'=\nonumber \\
 &= \sum_{j=1}^n \left(\sum_{i=1}^n b_{ij}{\Bbb Z}\right)\lambda_j'  &\subseteq  {\goth m}'.\nonumber 
\end{eqnarray}
Since both ${\goth m}'\subseteq {\goth m}$ and ${\goth m}\subseteq {\goth m}'$, we conclude
that ${\goth m}' = {\goth m}$. Lemma \ref{lm1} follows.
$\square$

\subsection{Projective invariance}
Recall that the measured foliations ${\cal F}$ and ${\cal F}'$ are 
{\it equivalent}, if there exists an automorphism $h\in Mod~(X)$,
which sends the leaves of the foliation ${\cal F}$ to the leaves of the foliation ${\cal F}'$.
This  equivalence deals with topological foliations, i.e.  projective
classes of measured foliations,  see  \cite{Thu1} for an explanation. 
\begin{lem}\label{lm2}
Let ${\cal F}, {\cal F}'$  be the equivalent measured foliations 
on a surface $X$. Then 
$$
Jac~({\cal F}')=\mu ~Jac~({\cal F}),
$$
where $\mu>0$ is a real number.  
\end{lem}
{\it Proof.} 
Let $h: X\to X$ be an automorphism of the surface $X$. Denote
by $h_*$ its action on $H_1(X, \Sing~({\cal F}); {\Bbb Z})$
and by $h^*$ on $H^1(X; {\Bbb R})$ connected  by the formula: 
$$
\int_{h_*(\gamma)}\omega=\int_{\gamma}h^*(\omega), \qquad\forall\gamma\in H_1(X, \Sing~({\cal F}); {\Bbb Z}), 
\qquad\forall\omega\in H^1(X; {\Bbb R}).
$$
Let $\omega,\omega'\in H^1(X; {\Bbb R})$ be the closed $1$-forms whose
trajectories define the foliations ${\cal F}$ and ${\cal F}'$, respectively.
Since ${\cal F}, {\cal F}'$ are equivalent measured foliations,
$$
\omega'= \mu ~h^*(\omega)
$$
for a $\mu>0$.

Let $Jac~({\cal F})={\Bbb Z}\lambda_1+\dots+{\Bbb Z}\lambda_n$ and 
$Jac~({\cal F}')={\Bbb Z}\lambda_1'+\dots+{\Bbb Z}\lambda_n'$. Then:
$$
\lambda_i'=\int_{\gamma_i}\omega'=\mu~\int_{\gamma_i}h^*(\omega)=
\mu~\int_{h_*(\gamma_i)}\omega, \qquad 1\le i\le n.
$$
By lemma \ref{lm1}, it holds:
$$
Jac~({\cal F})=\sum_{i=1}^n{\Bbb Z}\int_{\gamma_i}\omega=
\sum_{i=1}^n{\Bbb Z}\int_{h_*(\gamma_i)}\omega.
$$
Therefore:
$$
Jac~({\cal F}')=\sum_{i=1}^n{\Bbb Z}\int_{\gamma_i}\omega'=
\mu~\sum_{i=1}^n{\Bbb Z}\int_{h_*(\gamma_i)}\omega=\mu~Jac~({\cal F}).
$$
Lemma \ref{lm2} follows.
$\square$

\section{Equivalent foliations are stably isomorphic} 
Let ${\cal F}$ be a measured foliation on the surface $X$.
We introduce an AF $C^*$-algebra, ${\Bbb A}_{\cal F}$,  corresponding to  the foliation 
${\cal F}$ as explained in item G of introduction (for the foliation ${\cal F}_{\phi}$).
The goal of present section is the commutativity of the following  diagram:

\begin{picture}(300,110)(-80,-5)
\put(20,70){\vector(0,-1){35}}
\put(130,70){\vector(0,-1){35}}
\put(45,23){\vector(1,0){53}}
\put(45,83){\vector(1,0){53}}
\put(15,20){${\Bbb A}_{\cal F}$}
\put(128,20){${\Bbb A}_{{\cal F}'}$}
\put(17,80){${\cal F}$}
\put(125,80){${\cal F}'$}
\put(60,30){\sf stably}
\put(50,10){\sf isomorphic}

\put(54,90){\sf equivalent}
\end{picture}

\subsection{Modules and continued fractions}
The following lemma  is a simple 
property of the Jacobi-Perron fractions \cite{BE}. 
\begin{lem}\label{lm3}
Let ${\goth m}={\Bbb Z}\lambda_1+\dots+{\Bbb Z}\lambda_n$
and   ${\goth m}'={\Bbb Z}\lambda_1'+\dots+{\Bbb Z}\lambda_n'$
be two ${\Bbb Z}$-modules, such that ${\goth m}'=\mu {\goth m}$ for a $\mu>0$. 
Then the Jacobi-Perron continued fractions of the vectors $\lambda$ and $\lambda'$
coincide except,   possibly,  at a finite number of terms. 
\end{lem}
{\it Proof.}
Let ${\goth m}={\Bbb Z}\lambda_1+\dots+{\Bbb Z}\lambda_n$ and 
${\goth m}'={\Bbb Z}\lambda_1'+\dots+{\Bbb Z}\lambda_n'$. Since
${\goth m}'=\mu {\goth m}$, where $\mu$ is a positive real,
one gets the following identity of the ${\Bbb Z}$-modules:
$$
{\Bbb Z}\lambda_1'+\dots+{\Bbb Z}\lambda_n'={\Bbb Z}(\mu\lambda_1)+\dots+{\Bbb Z}(\mu\lambda_n).
$$
One can always assume that $\lambda_i$ and $\lambda_i'$ are positive reals.
For obvious reasons, there exists a basis $\{\lambda_1^{''},\dots,\lambda_n^{''}\}$
of the module ${\goth m}'$, such that:
$$
\left\{
\begin{array}{cc}
\lambda'' &= A(\mu\lambda) \nonumber\\
\lambda'' &= A'\lambda',
\end{array}
\right.
$$
where $A,A'\in GL^+_n({\Bbb Z})$ are the matrices, whose entries 
are non-negative integers.  In view of the Proposition 3 of  \cite{Bau1}:
$$
\left\{
\begin{array}{cc}
A &=  \left(\matrix{0 & 1\cr I & b_1}\right)\dots
\left(\matrix{0 & 1\cr I & b_k}\right)\nonumber\\
A' &= \left(\matrix{0 & 1\cr I & b_1'}\right)\dots
\left(\matrix{0 & 1\cr I & b_l'}\right),
\end{array}
\right.
$$
where $b_i, b_i'$ are non-negative integer vectors.
Since the  (Jacobi-Perron) continued fraction for the vectors
$\lambda$ and $\mu\lambda$ coincide for any $\mu>0$ \cite{BE},
we conclude that: 
$$
\left\{
\begin{array}{cc}
\left(\matrix{1\cr \theta}\right)
 &=  \left(\matrix{0 & 1\cr I & b_1}\right)\dots
\left(\matrix{0 & 1\cr I & b_k}\right)
\left(\matrix{0 & 1\cr I & a_1}\right)
\left(\matrix{0 & 1\cr I & a_2}\right)\dots
\left(\matrix{0\cr {\Bbb I}}\right)
\nonumber\\
\left(\matrix{1\cr \theta'}\right)
 &= \left(\matrix{0 & 1\cr I & b_1'}\right)\dots
\left(\matrix{0 & 1\cr I & b_l'}\right)
\left(\matrix{0 & 1\cr I & a_1}\right)
\left(\matrix{0 & 1\cr I & a_2}\right)\dots
\left(\matrix{0\cr {\Bbb I}}\right),
\end{array}
\right.
$$
where 
$$
\left(\matrix{1\cr \theta''}\right)=
\lim_{i\to\infty} \left(\matrix{0 & 1\cr I & a_1}\right)\dots
\left(\matrix{0 & 1\cr I & a_i}\right)
\left(\matrix{0\cr {\Bbb I}}\right). 
$$
In other words, the continued fractions of the vectors $\lambda$ and $\lambda'$
coincide  except at a finite number of terms.
$\square$

\subsection{Main lemma}
\begin{lem}\label{lm4}
Let ${\cal F}$ and  $ {\cal F}'$ be  equivalent measured foliations on a surface $X$.
Then the AF $C^*$-algebras ${\Bbb A}_{\cal F}$ and  ${\Bbb A}_{{\cal F}'}$
are stably isomorphic.  
\end{lem}
{\it Proof.}
Notice that lemma \ref{lm2} implies that
 equivalent measured foliations ${\cal F}, {\cal F}'$ have
 proportional jacobians, i.e. ${\goth m}'=\mu {\goth m}$
for a $\mu>0$.  On the other hand, by lemma \ref{lm3}
the continued fraction expansion of the basis vectors
of the proportional jacobians must coincide, except
a finite number of terms. Thus, the AF $C^*$-algebras
${\Bbb A}_{\cal F}$ and ${\Bbb A}_{{\cal F}'}$ 
are given by  the Bratteli diagrams, which are identical,
except a finite part of the diagram.  
It is well  known (\cite{E}, Theorem 2.3)  that the AF $C^*$-algebras,  
which have such a  property, are stably isomorphic.      
$\square$

\section{Proofs} 
\subsection{Proof of theorem 1} 
Let $\phi\in Mod~(X)$ be a pseudo-Anosov automorphism of the surface $X$.
Denote by ${\cal F}_{\phi}$ the invariant foliation of $\phi$.
By definition of such a foliation, $\phi({\cal F}_{\phi})=\lambda_{\phi}{\cal F}_{\phi}$,
where $\lambda_{\phi}>1$ is the dilatation of $\phi$.

Consider the jacobian $Jac~({\cal F}_{\phi})={\goth m}_{\phi}$
of  foliation ${\cal F}_{\phi}$. 
Since ${\cal F}_{\phi}$ is an invariant foliation of the pseudo-Anosov automorphism $\phi$, 
one gets the following equality  of the ${\Bbb Z}$-modules:
\begin{equation}\label{eq1}
{\goth m}_{\phi}=\lambda_{\phi}{\goth m}_{\phi}, \qquad \lambda_{\phi}\ne\pm 1.
\end{equation} 
Let $\{v^{(1)},\dots,v^{(n)}\}$ be a basis in  module  ${\goth m}_{\phi}$,
such that $v^{(i)}>0$. In view of (\ref{eq1}), one obtains the following
system of linear equations:
\begin{equation}\label{eq2}
\left\{
\begin{array}{ccc}
\lambda_{\phi}v^{(1)}  &=& a_{11}v^{(1)}+a_{12}v^{(2)}+\dots+a_{1n}v^{(n)}\\
\lambda_{\phi}v^{(2)}  &=& a_{21}v^{(1)}+a_{22}v^{(2)}+\dots+a_{2n}v^{(n)}\\
\vdots && \\
\lambda_{\phi}v^{(n)}  &=& a_{n1}v^{(1)}+a_{n2}v^{(2)}+\dots+a_{nn}v^{(n)},
\end{array}
\right.
\end{equation}
where $a_{ij}\in {\Bbb Z}$. The matrix $A=(a_{ij})$ is invertible. Indeed,
since  foliation ${\cal F}_{\phi}$ is minimal, real numbers  $v^{(1)},\dots,v^{(n)}$ 
are linearly independent over ${\Bbb Q}$. So do numbers $\lambda_{\phi}v^{(1)},\dots,\lambda_{\phi}v^{(n)}$,
which therefore can be taken for a basis of the module ${\goth m}_{\phi}$. 
Thus, there exists an integer matrix $B=(b_{ij})$, such that $v^{(j)}=\sum_{i,j}w^{(i)}$,
where $w^{(i)}=\lambda_{\phi}v^{(i)}$. Clearly, $B$ is an inverse to  matrix $A$.
Therefore, $A\in GL_n({\Bbb Z})$.

Moreover, without loss of the generality one can assume that $a_{ij}\ge0$. 
Indeed, if it is not yet the case, consider the conjugacy class $[A]$
of the matrix $A$. Since $v^{(i)}>0$,  there exists a matrix $A^+\in [A]$
whose entries are non-negative integers. One has to replace 
 basis $v=(v^{(1)},\dots,v^{(n)})$ in the module ${\goth m}_{\phi}$ 
by a  basis  $Tv$, where $A^+=TAT^{-1}$. It will be further assumed that
 $A=A^+$. 
\begin{lem}\label{lm5}
Vector $(v^{(1)}, \dots, v^{(n)})$ is the  limit of a periodic  
Jacobi-Perron continued fraction. 
\end{lem}
{\it Proof.} It follows from the discussion above,  that there exists
a non-negative integer matrix $A$,  such that $Av=\lambda_{\phi}v$. 
In view of \cite{Bau1},  Proposition 3,  matrix $A$ admits 
a unique factorization:
\begin{equation}\label{eq3}
A=
\left(\matrix{0 & 1\cr I & b_1}\right)\dots
\left(\matrix{0 & 1\cr I & b_k}\right),
\end{equation}
where $b_i=(b^{(i)}_1,\dots, b^{(i)}_{n})^T$ are vectors of the  non-negative integers.
Let us consider the  periodic Jacobi-Perron continued fraction:
\begin{equation}\label{eq4}
Per
~\overline{
 \left(\matrix{0 & 1\cr I & b_1}\right)\dots
\left(\matrix{0 & 1\cr I & b_k}\right)
}
\left(\matrix{0\cr {\Bbb I}}\right).
\end{equation}
According to \cite{Per1}, {\bf Satz XII},  the above  fraction converges to vector $w=(w^{(1)},\dots, w^{(n)})$,
such that $w$ satisfies equation $(B_1B_2\dots B_k)w=Aw=\lambda_{\phi}w$.
In view of  equation $Av=\lambda_{\phi}v$, we conclude that  vectors $v$ and $w$
are collinear.  Therefore, the Jacobi-Perron continued fractions of $v$ and $w$
must coincide.  
$\square$

\bigskip\noindent
It is now straightforward to prove, that the AF $C^*$-algebra attached to  foliation
${\cal F}_{\phi}$ is stationary.  Indeed, by lemma \ref{lm5}, 
the vector of periods $v^{(i)}=\int_{\gamma_i}\omega$ unfolds  into a periodic Jacobi-Perron 
continued fraction. By the definition, the Bratteli diagram of the AF $C^*$-algebra  ${\Bbb A}_{\phi}$
is periodic as well.  In other words, the AF $C^*$-algebra  ${\Bbb A}_{\phi}$ is stationary. 
$\square$

\subsection{Proof of theorem 2} 
(i) Let us prove the first statement. For the sake of completeness,
let us give a proof of the following (well-known) lemma. 
\begin{lem}\label{lm6}
Let $\phi$ and $\phi'$ be  conjugate pseudo-Anosov automorphisms
of a surface $X$. Then their invariant foliations ${\cal F}_{\phi}$
and ${\cal F}_{\phi'}$ are equivalent  as measured foliations.  
\end{lem}
{\it Proof.} Let $\phi,\phi'\in Mod~(X)$ be conjugate, i.e 
$\phi'=h\circ\phi\circ h^{-1}$ for an automorphism $h\in Mod~(X)$.
Since $\phi$ is the pseudo-Anosov automorphism, there exists  a measured foliation
${\cal F}_{\phi}$,  such that $\phi({\cal F}_{\phi})=\lambda_{\phi}{\cal F}_{\phi}$.
Let us evaluate the automorphism $\phi'$ on the foliation $h({\cal F}_{\phi})$:
\begin{eqnarray}\label{eq5}
\phi'(h({\cal F}_{\phi}))  &= h\phi h^{-1}(h({\cal F}_{\phi})) &= 
h\phi({\cal F}_{\phi})=\nonumber \\
 &= h \lambda_{\phi} {\cal F}_{\phi}  &= \lambda_{\phi} (h({\cal F}_{\phi})). 
\end{eqnarray}
Thus, ${\cal F}_{\phi'}=h({\cal F}_{\phi})$ is the invariant foliation for the 
pseudo-Anosov automorphism $\phi'$ and ${\cal F}_{\phi}, {\cal F}_{\phi'}$
are equivalent foliations. Note also that the pseudo-Anosov automorphism $\phi'$ has 
the same dilatation as the automorphism $\phi$.  
$\square$

\bigskip\noindent
Suppose that $\phi$ and $\phi'$ are  conjugate pseudo-Anosov 
automorphisms.  Functor $F$ acts by the formulas $\phi\mapsto {\Bbb A}_{\phi}$
and $\phi'\mapsto {\Bbb A}_{\phi'}$, where ${\Bbb A}_{\phi}, {\Bbb A}_{\phi'}$
are the AF $C^*$-algebras corresponding to   invariant foliations ${\cal F}_{\phi}, {\cal F}_{\phi'}$. 
In view of lemma \ref{lm6}, ${\cal F}_{\phi}$ and ${\cal F}_{\phi'}$ are
 equivalent measured foliations.  Then, by lemma \ref{lm4}, the AF $C^*$-algebras ${\Bbb A}_{\phi}$
and ${\Bbb A}_{\phi'}$ are stably isomorphic AF $C^*$-algebras.  Item (i)
follows.

\bigskip
(ii) Let us prove the second statement. We start with an elementary observation.
Let $\phi\in Mod~(X)$ be a pseudo-Anosov automorphism. Then there exists a unique
measured foliation, ${\cal F}_{\phi}$, such that $\phi({\cal F}_{\phi})=\lambda_{\phi}{\cal F}_{\phi}$,
where $\lambda_{\phi}>1$ is an algebraic integer. Let us evaluate  automorphism
$\phi^2\in Mod~(X)$ on the foliation ${\cal F}_{\phi}$: 
\begin{eqnarray}\label{eq6}
\phi^2({\cal F}_{\phi}) &= \phi (\phi({\cal F}_{\phi})) &= 
\phi(\lambda_{\phi} {\cal F}_{\phi})=\nonumber \\
= \lambda_{\phi} \phi({\cal F}_{\phi}) &= \lambda_{\phi}^2{\cal F}_{\phi}  &= 
\lambda_{\phi^2}{\cal F}_{\phi}, 
\end{eqnarray}
where $\lambda_{\phi^2}:= \lambda_{\phi}^2$. Thus,  foliation ${\cal F}_{\phi}$
is an invariant foliation for the automorphism $\phi^2$ as well. By induction,
one concludes that ${\cal F}_{\phi}$ is an invariant foliation of  the 
automorphism $\phi^n$ for any $n\ge 1$.

Even more is true. Suppose that $\psi\in Mod~(X)$ is a pseudo-Anosov
automorphism, such that $\psi^m=\phi^n$ for some $m\ge 1$ and $\psi\ne\phi$.
Then ${\cal F}_{\phi}$ is an invariant foliation for the automorphism 
$\psi$. Indeed, ${\cal F}_{\phi}$ is  invariant foliation of  the 
automorphism $\psi^m$. If there exists ${\cal F}'\ne {\cal F}_{\phi}$
such that the foliation  ${\cal F}'$ is an invariant foliation  of $\psi$, then 
the foliation ${\cal F}'$ is
also an invariant foliation  of the pseudo-Anosov automorphism $\psi^m$. 
Thus, by the uniqueness, 
${\cal F}'={\cal F}_{\phi}$. We have just proved the following lemma. 
\begin{lem}\label{lm7}
Let $\phi$ be the pseudo-Anosov automorphism of a surface $X$.
Denote by $[\phi]$ a set of the pseudo-Anosov automorphisms $\psi$
of $X$,  such that $\psi^m=\phi^n$ for some positive integers
$m$ and $n$. Then the pseudo-Anosov foliation ${\cal F}_{\phi}$
is an invariant foliation for every pseudo-Anosov automorphism $\psi\in [\phi]$. 
\end{lem}

\medskip\noindent
In view of lemma \ref{lm7}, one arrives at the following
identities among the AF $C^*$-algebras:
\begin{equation}\label{eq7}
{\Bbb A}_{\phi}={\Bbb A}_{\phi^2}=\dots={\Bbb A}_{\phi^n}=
{\Bbb A}_{\psi^m}=\dots={\Bbb A}_{\psi^2}={\Bbb A}_{\psi}.
\end{equation}
Thus,   functor $F$ is not an 
injective functor: the preimage, $Ker~F$, of  algbera
${\Bbb A}_{\phi}$ consists of a countable set of the pseudo-Anosov
automorphisms $\psi\in [\phi]$, commensurable with the automorphism
$\phi$.

\medskip
Theorem \ref{thm2} is proved.
$\square$

\subsection{Proof of corollary 1}
(i) It follows from theorem \ref{thm1}, that ${\Bbb A}_{\phi}$
is a stationary AF $C^*$-algebra. An arithmetic invariant of the stable
isomorphism classes of the stationary AF $C^*$-algebras has been found 
by D.~Handelman in \cite{Han1}. Summing up his results, the invariant
is as follows. 

Let $A\in GL_n({\Bbb Z})$ be a matrix with the strictly positive
entries, such that $A$ is  equal to the minimal period of the Bratteli diagram
of the stationary AF $C^*$-algebra. (In case the matrix $A$ has  zero entries,
it is necessary to take a proper minimal power of the matrix $A$.) By the Perron-Frobenius
theory,  matrix $A$ has a real eigenvalue $\lambda_A>1$, which exceeds
the absolute values of other roots of the characteristic polynomial of $A$.
Note that $\lambda_A$ is an invertible algebraic integer (the unit).  Consider
the  real algebraic number field $K={\Bbb Q}(\lambda_A)$ obtained as 
an extension of the field of the rational numbers by the algebraic 
number $\lambda_A$. Let $(v^{(1)}_A,\dots,v^{(n)}_A)$ be the eigenvector
corresponding to the eigenvalue $\lambda_A$. One can  normalize the eigenvector so 
that $v^{(i)}_A\in K$.

The   departure point  of Handelman's invariant is the  ${\Bbb Z}$-module
${\goth m}={\Bbb Z}v^{(1)}_A+\dots+{\Bbb Z}v^{(n)}_A$. The module ${\goth m}$
brings in two new arithmetic objects: (i) the ring $\Lambda$ of the endomorphisms of ${\goth m}$  
and (ii) an ideal $I$ in the ring $\Lambda$,  such that $I={\goth m}$  after a scaling (\cite{BS}, Lemma 1, p.88). 
The ring $\Lambda$ is an order in the algebraic number field $K$ and therefore one can talk about the ideal
classes in $\Lambda$. The ideal class of $I$ is denoted by $[I]$. 
Omitting the embedding question for the field $K$,  the triple $(\Lambda, [I], K)$ is an invariant of
the  stable isomorphism class of the  stationary AF $C^*$-algebra ${\Bbb A}_{\phi}$  (\S 5 of \cite{Han1}).  
Item (i) follows.

\bigskip
(ii)  Numerical invariants of the stable isomorphism
classes of the stationary AF $C^*$-algebras can be  derived from the  triple $(\Lambda, [I], K)$.
These invariants are the rational integers -- called the determinant and signature -- can be  obtained as 
follows.

Let ${\goth m}, {\goth m}'$ be the full ${\Bbb Z}$-modules in an algebraic
number field $K$. It follows from (i), that if ${\goth m}\ne {\goth m}'$
are distinct as the ${\Bbb Z}$-modules,
then the corresponding AF $C^*$-algebras cannot be stably isomorphic. 
We wish to find the numerical invariants, which discern the case ${\goth m}\ne {\goth m}'$. 
It is assumed that a ${\Bbb Z}$-module is given by the  set of generators
$\{\lambda_1,\dots,\lambda_n\}$. Therefore, the  problem can be formulated as follows: find a number
attached to the set of generators $\{\lambda_1,\dots,\lambda_n\}$,  which does not change   
on the set of generators $\{\lambda_1',\dots,\lambda_n'\}$ of  the same {\Bbb Z}-module.

One such invariant is associated with the trace function on the algebraic number field $K$. 
Recall that $Tr: K\to {\Bbb Q}$ is a linear function on  $K$ such
that $Tr ~(\alpha+\beta)=\tr~(\alpha)+ \tr~(\beta)$ and
$\tr~(a\alpha)=a ~\tr~(\alpha)$ for $\forall\alpha,\beta\in K$ and
$\forall a\in {\Bbb Q}$.

Let ${\goth m}$ be a full ${\Bbb Z}$-module in the field $K$.   
The trace function defines a symmetric bilinear form
$q(x,y): {\goth m}\times {\goth m}\to {\Bbb Q}$ by 
the formula:
\begin{equation}\label{eq8}
(x,y)\longmapsto \tr~(xy), \qquad \forall x,y\in {\goth m}.
\end{equation}
The form $q(x,y)$ depends on the basis $\{\lambda_1,\dots,\lambda_n\}$
in the module ${\goth m}$:
\begin{equation}\label{eq16}
q(x,y)=\sum_{j=1}^n\sum_{i=1}^na_{ij}x_iy_j, \qquad\hbox{where} \quad a_{ij}=\tr~(\lambda_i\lambda_j).
\end{equation}
However, the general theory of the bilinear forms (over the fields ${\Bbb Q}, {\Bbb R}, {\Bbb C}$ or the ring
of rational integers ${\Bbb Z}$) 
tells us that certain numerical quantities will not depend on the choice of such a basis.

Namely, one such invariant is as follows. Consider a symmetric matrix $A$
corresponding to the bilinear form $q(x,y)$:
\begin{equation}\label{eq9}
A=\left(
\matrix{a_{11} & a_{12} & \dots & a_{1n}\cr
        a_{12} & a_{22} & \dots & a_{2n}\cr
        \vdots  &        &       & \vdots\cr
        a_{1n} & a_{2n} & \dots & a_{nn}}
\right).
\end{equation}
It is known that the matrix $A$, written in a new basis, will take the
form $A'=U^TAU$, where ~$U\in GL_n({\Bbb Z})$. 
Then $\dete~(A')=\dete~(U^TAU)=$
\linebreak
$=\dete~(U^T) \dete~(A) \dete~(U) =\dete~(A)$.
Therefore, the rational integer number:
\begin{equation}\label{eq10}
\Delta= \dete~(\tr~(\lambda_i\lambda_j)),
\end{equation}
called a {\it determinant} of the bilinear form  $q(x,y)$, does not depend on the
choice of the basis $\{\lambda_1,\dots,\lambda_n\}$ in the
module ${\goth m}$.  We conclude that the determinant 
$\Delta$ discerns
\footnote{Note that if $\Delta=\Delta'$ for the modules ${\goth m}, {\goth m}'$,
one cannot conclude that ${\goth m}={\goth m}'$. The problem of equivalence
of the  symmetric bilinear forms over ${\Bbb Q}$ (i.e. the existence of a linear
substitution over ${\Bbb Q}$,  which transforms one form to the other),
is a fundamental question of number theory. The Minkowski-Hasse theorem
says that two such forms are equivalent if and only if they are
equivalent over the field ${\Bbb Q}_p$ for every prime number $p$
and over the field ${\Bbb R}$. Clearly, the resulting $p$-adic quantities
will give new invariants of the stable isomorphism classes of the AF $C^*$-algebras.
The question is much similar to the Minkowski units attached to knots, see e.g.
Reidemeister \cite{R}. We will not pursue this topic here and refer the reader
to the problem part of present article.}
the modules ${\goth m}\ne {\goth m}'$.

\medskip
Finally, recall that the form $q(x,y)$ can be brought by an integer linear
transformation  to the diagonal form:
\begin{equation}\label{eq11}
a_1x_1^2+a_2x_2^2+\dots+a_nx_n^2,
\end{equation}
where $a_i\in {\Bbb Z}\setminus \{0\}$. We let $a_i^+$ be the positive and $a_i^-$ the 
negative entries in the diagonal form.   In view of the law of inertia
for the bilinear forms, the integer number
$\Sigma = (\# a_i^+) - (\# a_i^-)$,
called  a {\it signature}, does not depend
on a  particular choice of the basis in the module ${\goth m}$.
Thus, $\Sigma$ discerns the modules ${\goth m}\ne {\goth m}'$.
Corollary \ref{cor1} follows.
$\square$

\section{Examples, open problems and conjectures}
In the present section we shall calculate  invariants $\Delta$ and $\Sigma$
for the Anosov automorphisms of the two-dimensional torus.  
Examples of two non-conjugate Anosov automorphisms with
the same Alexander polynomial, but different determinants  $\Delta$ are 
constructed.  Recall that isotopy classes of the orientation-preserving diffeomorphisms
of the torus $T^2$ are bijective with the $2\times 2$ matrices with integer entries
and determinant $+1$, i.e. $Mod~ (T^2)\cong SL(2,{\Bbb Z})$. Under the identification,
the non-periodic automorphisms correspond to the matrices $A\in SL(2,{\Bbb Z})$ 
with $|\tr~A|>2$.

\subsection{Full modules and orders  in the quadratic field}
Let $K={\Bbb Q}(\sqrt{d})$ be a quadratic extension of the field of rational numbers ${\Bbb Q}$.
Further we suppose that $d$ is a positive square free integer. Let
\begin{equation}
\omega=\cases{{1+\sqrt{d}\over 2} & if $d\equiv 1 ~mod~4$,\cr
               \sqrt{d} & if $d\equiv 2,3 ~mod~4$.}
\end{equation}
\begin{prp}\label{prp1}
Let $f$ be a positive integer. 
Every order  in $K$ has form  
$\Lambda_f={\Bbb Z} +(f\omega){\Bbb Z}$, where 
$f$ is the conductor of $\Lambda_f$.
\end{prp}
{\it Proof.} See \cite{BS} pp. 130-132.
$\square$

\medskip\noindent 
The proposition \ref{prp1} allows to classify the similarity classes of the full modules 
in the field $K$. Indeed, there exists a finite number of ${\goth m}_f^{(1)},\dots,
{\goth m}_f^{(s)}$ of the non-similar full modules in the field  $K$, whose coefficient
ring is the order $\Lambda_f$, cf Theorem 3, Ch 2.7 of \cite{BS}. Thus, proposition \ref{prp1}
gives a finite-to-one classification of the similarity classes of full modules in the field $K$.

\subsection{Numerical invariants of the Anosov automorphisms}
Let $\Lambda_f$ be an order  in $K$ with
the  conductor $f$. Under the addition operation, the order  $\Lambda_f$ is a full module, 
which we denote by ${\goth m}_f$.  
Let us evaluate  the invariants $q(x,y)$, $\Delta$  and $\Sigma$ 
on the module  ${\goth m}_f$. To calculate  $(a_{ij})=\tr~(\lambda_i\lambda_j)$,  we let 
$\lambda_1=1,\lambda_2=f\omega$. Then:
\begin{eqnarray}
a_{11} &=& 2, \quad a_{12}=a_{21}= f, \quad a_{22}= {1\over 2} f^2(d+1)\quad \hbox{if} \quad d\equiv 1 ~mod~4\nonumber\\
a_{11} &=&  2, \quad a_{12}=a_{21}= 0, \quad a_{22}= 2f^2d \quad \hbox{if} \quad d\equiv 2,3 ~mod~4, 
\end{eqnarray}
and 
\begin{eqnarray}\label{eq15+}
q(x,y) &=& 2x^2 +2f xy +{1\over 2}f^2(d+1)y^2\quad \hbox{if} \quad d\equiv 1 ~mod~4\nonumber\\
q(x,y) &=&  2x^2+2f^2dy^2\quad  \hbox{if} \quad d\equiv 2,3 ~mod~4.
\end{eqnarray}
Therefore
\begin{equation}\label{eq16+}
\Delta=\cases{f^2d  & if $d\equiv 1 ~mod~4$,\cr
               4f^2d & if $d\equiv 2,3 ~mod~4$,}
\end{equation}
and $\Sigma=+2$ in the both cases, where $\Sigma=\# (positive) -\# (negative)$ 
entries in the diagonal normal form of $q(x,y)$.

\subsection{Examples}
Let us consider some numerical examples,  which illustrate 
advantages  of our  invariants in comparison to the classical Alexander polynomials.  
\begin{exm}\label{ex1}
\textnormal{
Denote by $M_A$ and $M_B$ the hyperbolic 3-dimensional manifolds obtained as a torus
bundle over the circle with the monodromies
\begin{equation}
A=\left(\matrix{5 & 2\cr 2 & 1}\right)\qquad  \hbox{and}
\qquad B=\left(\matrix{5 & 1\cr 4 & 1}\right), 
\end{equation}
respectively. The Alexander polynomials of $M_A$ and $M_B$ are identical   
$\Delta_A(t)=\Delta_B(t)= t^2-6t+1$. However, the manifolds $M_A$ and $M_B$ are 
{\it not} homotopy equivalent. 
Indeed, the Perron-Frobenius eigenvector of matrix $A$ is  $v_A=(1, \sqrt{2}-1)$
while of the matrix $B$ is $v_B=(1, 2\sqrt{2}-2)$. The  bilinear forms for the modules 
${\goth m}_A={\Bbb Z}+(\sqrt{2}-1){\Bbb Z}$ and 
${\goth m}_B={\Bbb Z}+(2\sqrt{2}-2){\Bbb Z}$ can be written as
\begin{equation}
q_A(x,y)= 2x^2-4xy+6y^2,\qquad q_B(x,y)=2x^2-8xy+24y^2,
\end{equation}
respectively.  The  modules ${\goth m}_A, {\goth m}_B$ are not similar in the number field 
$K={\Bbb Q}(\sqrt{2})$, since   their  determinants $\Delta({\goth m}_A)=8$ and 
$\Delta({\goth m}_B)=32$ are not equal. Therefore,    matrices $A$ and $B$ are not conjugate
\footnote{The reader may verify this fact using the method of periods, which dates back
to Gauss. First we have to find the fixed points  $Ax=x$ and $Bx=x$, which gives us $x_A=1+\sqrt{2}$ and
$x_B={1+\sqrt{2}\over 2}$, respectively. Then one unfolds the fixed points into a periodic continued fraction,
which gives us $x_A=[2,2,2,\dots]$ and $x_B=[1,4,1,4,\dots]$. Since the period $(2)$ of $x_A$
differs from the period $(1,4)$ of $B$,  the matrices $A$ and $B$ belong to different
conjugacy classes in $SL(2, {\Bbb Z})$.}
in the group  $SL(2,{\Bbb Z})$. Note that the class number $h_K=1$ for the field $K$.     
}
\end{exm}

\begin{exm}\label{ex2}
{\bf (\cite{Han2}, p.12)}
\textnormal{
Let  $M_A$ and $M_B$ be  3-dimensional manifolds corresponding to
matrices
\begin{equation}
A=\left(\matrix{4 & 3\cr 5 & 4}\right)\qquad  \hbox{and}
\qquad B=\left(\matrix{4 & 15\cr 1 & 4}\right), 
\end{equation}
respectively.  The Alexander polynomials of $M_A$ and $M_B$ are identical   
$\Delta_A(t)=\Delta_B(t)= t^2-8t+1$.   Yet  the manifolds $M_A$ and $M_B$ are 
not homotopy equivalent. 
Indeed, the Perron-Frobenius eigenvector of matrix $A$ is  $v_A=(1,  {1\over 3}\sqrt{15})$
while of the matrix $B$ is $v_B=(1, {1\over 15}\sqrt{15})$.  
The corresponding  modules are  ${\goth m}_A={\Bbb Z}+( {1\over 3}\sqrt{15}){\Bbb Z}$ and 
${\goth m}_B={\Bbb Z}+( {1\over 15}\sqrt{15}){\Bbb Z}$;  note that  $d=15\equiv 3~mod~4$
in both cases, but the corresponding conductors are  $f_A=3$ and $f_B=15$. 
Using formulas (\ref{eq15+}) one finds
\begin{equation}
q_A(x,y)= 2x^2+18y^2,\qquad q_B(x,y)=2x^2 +450y^2,
\end{equation}
respectively.  The  modules ${\goth m}_A, {\goth m}_B$ are not similar in the number field 
$K={\Bbb Q}(\sqrt{15})$, since  formulas (\ref{eq16+}) imply that  their  determinants $\Delta({\goth m}_A)=36$ and 
$\Delta({\goth m}_B)=900$ are not equal. Therefore,  matrices $A$ and $B$ are not conjugate
in the group  $SL(2,{\Bbb Z})$.     
}
\end{exm}

\begin{exm}\label{ex3}
{\bf (\cite{Han2}, p.12})
\textnormal{
Let $a, b$ be a pair of positive integers satisfying the Pell equation
$a^2-8b^2=1$;  the latter has infinitely many  solutions, e.g. $a=3, b=1$,
{\it etc.}   Denote by   $M_A$ and $M_B$ the  3-dimensional manifolds corresponding to
matrices
\begin{equation}
A=\left(\matrix{a & 4b\cr 2b & a}\right)\qquad  \hbox{and}
\qquad B=\left(\matrix{a & 8b\cr b & a}\right), 
\end{equation}
respectively. 
The Alexander polynomials of $M_A$ and $M_B$ are identical   
$\Delta_A(t)=\Delta_B(t)= t^2-2a t+1$.   Yet  manifolds $M_A$ and $M_B$ are 
not homotopy equivalent. 
Indeed, the Perron-Frobenius eigenvector of matrix $A$ is  $v_A=(1,  {1\over 4b}\sqrt{a^2-1})$
while of the matrix $B$ is $v_B=(1, {1\over 8b}\sqrt{a^2-1})$.  
The corresponding  modules are  ${\goth m}_A={\Bbb Z}+( {1\over 4b}\sqrt{a^2-1}){\Bbb Z}$ and 
${\goth m}_B={\Bbb Z}+( {1\over 8b}\sqrt{a^2-1}){\Bbb Z}$.  It is easy to see  that  the discriminant $d=a^2-1\equiv 3~mod~4$
for all  $a\ge 2$.  Indeed, $d=(a-1)(a+1)$ and,   therefore,  integer   $a\not\equiv 1; 3 ~mod~4$;  
hence $a\equiv 2~mod~4$ so that  $a-1\equiv 1~mod ~4$
and $a+1\equiv 3 ~mod  ~4$ and, thus,  $d=a^2-1\equiv 3~mod~4$.    
 Therefore the  corresponding conductors are  $f_A=4b$ and $f_B=8b$, 
and
\begin{equation}
q_A(x,y)= 2x^2+32b^2(a^2-1)y^2,\quad q_B(x,y)=2x^2 +128b^2(a^2-1)y^2,
\end{equation}
respectively.  The  modules ${\goth m}_A, {\goth m}_B$ are not similar in the number field 
$K={\Bbb Q}(\sqrt{a^2-1})$, since   their  determinants $\Delta({\goth m}_A)=64~b^2(a^2-1)$ and 
$\Delta({\goth m}_B)=256~b^2(a^2-1)$ are not equal.  Therefore,  matrices $A$ and $B$ are not conjugate
in the group  $SL(2,{\Bbb Z})$. 
}
\end{exm}


\subsection{Open problems and conjectures}
This  section is reserved for some  questions and conjectures, 
which arise in  connection with the invariants 
$(\Lambda, [I], K), q(x,y), \Delta$ and $\Sigma$.

\bigskip
{\it 1. $P$-adic invariants of the pseudo-Anosov automorphisms}

\medskip\noindent
{\bf A.} Let $\phi\in Mod~(X)$ be pseudo-Anosov automorphism
of a surface $X$. If $\lambda_{\phi}$ is the dilatation of
$\phi$, then one can consider a ${\Bbb Z}$-module 
${\goth m}={\Bbb Z}v^{(1)}+\dots+{\Bbb Z}v^{(n)}$ in the
number field $K={\Bbb Q}(\lambda_{\phi})$ generated by the
normalized eigenvector $(v^{(1)},\dots,v^{(n)})$ corresponding
to the eigenvalue $\lambda_{\phi}$. The trace function on the number field $K$
gives rise to a symmetric bilinear form $q(x,y)$ on the 
module ${\goth m}$. The form is defined over the field ${\Bbb Q}$. 
It has been shown that a pseudo-Anosov automorphism $\phi'$,
conjugate to $\phi$, yields a form $q'(x,y)$, equivalent to
$q(x,y)$, i.e. $q(x,y)$ can be transformed to $q'(x,y)$ by
an invertible linear substitution with the coefficients in ${\Bbb Z}$.

\medskip\noindent
{\bf B.} Recall that two rational bilinear forms
$q(x,y)$ and  $q'(x,y)$ are equivalent whenever 
the following conditions are met:

\medskip
(i) $\Delta=\Delta'$, where $\Delta$ is the determinant of the form;

\smallskip
(ii) for each prime number $p$ (including $p=\infty$) certain $p$-adic
equation between the coefficients of forms $q, q'$ must be satisfied,
see e.g. \cite{BS}, Ch.1, \S 7.5. (In fact, only a {\it finite} 
number of such equations have to be verified.)

\medskip\noindent
Condition (i) has been already used to discern between the 
conjugacy classes of the  pseudo-Anosov automorphisms. One can use
condition (ii) to discern between the pseudo-Anosov automorphisms
with $\Delta=\Delta'$. The following question can be posed: 
{\it Find the $p$-adic invariants of the pseudo-Anosov automorphisms.}

\bigskip
{\it 2. Signature of the  pseudo-Anosov automorphism}

\medskip\noindent
The signature is an important and well-known invariant connected to
the chirality and knotting number of knots and links \cite{R}.
It will be interesting to find a geometric interpretation
of  the signature $\Sigma$ for  the pseudo-Anosov automorphisms. 
One can ask the following  question: 
{\it Find a geometric meaning of the invariant $\Sigma$.}

\bigskip
{\it 3.  Number of conjugacy classes of the  pseudo-Anosov automorphisms
with the  same dilatation}

\medskip\noindent
The dilatation $\lambda_{\phi}$ is an invariant of the conjugacy class of
the pseudo-Anosov automorphism $\phi\in Mod~(X)$. On the other hand, 
it is known that there exist non-conjugate pseudo-Anosov's with 
the same dilatation and the number of such classes is finite \cite{Thu1}. 
It is natural to expect that the invariants of operator algebras
can be used to evaluate the number. We conclude with the following
conjecture.  
\begin{con}\label{cn1}
Let $(\Lambda, [I], K)$ be the triple corresponding to a pseudo-Anosov
automorphism $\phi\in Mod~(X)$. Then the number of the conjugacy classes
of the pseudo-Anosov automorphisms with the dilatation $\lambda_{\phi}$
is equal to the class number $h_{\Lambda}=|\Lambda/[I]|$ 
of the integral order $\Lambda$. 
\end{con}

\bigskip\noindent
{\sf Acknowledgment.} 
I thank the referee for helpful comments and Daniel Silver and Susan 
Williams for their interest and hospitality.

    

\vskip1cm

\textsc{The Fields Institute for Mathematical Sciences, Toronto, ON, Canada,  
E-mail:} {\sf igor.v.nikolaev@gmail.com}

\smallskip
{\it Current address: 616-315 Holmwood Ave., Ottawa, ON, Canada, K1S 2R2}

\end{document}